\documentstyle{amsppt}
\baselineskip=2\baselineskip
\def\annl#1{{}_{\g A}\!\operatorname{ann}#1}
\def\annr#1{\text{\rm ann}_{\g A}\,#1}
\def\ann#1{\operatorname{ann}#1}
\def\card#1{\operatorname{card}(#1)}

\def\DC#1{\pmb\Delta(#1)}
\def\norm#1{\Vert #1\Vert}
\def\sbst{\subseteq}

\def\wt{\widetilde}

\define\ha#1{\bold h_{\g A}(\g A,#1)}
\define\ah#1{{}_{\g A}\bold h(\g A,#1)}
\define\aha#1{{}_{\g A}\bold h_{\g A}(\g A,#1)}
\define\g{\goth}
\define\A#1{\Cal A(#1)}
\define\F#1{\Cal F(#1)}
\define\N#1{\Cal N(#1)}
\define\I#1{\Cal I(#1)}
\define\Si#1{\Cal S(#1)}
\define\T#1{#1\widehat\otimes#1^*}
\define\Bn#1{\Cal B(#1)}
\define\Aa#1#2{\Cal A(#1,#2)}
\define\Ff#1#2{\Cal F(#1,#2)}
\define\nul{\{0\}}
\define\flp#1{\flushpar{\rm(#1)}}
\define\Ii#1#2{\Cal I(#1,#2)}
\define\Nn#1#2{\Cal N(#1,#2)}
\define\Bb#1#2{\Cal B(#1,#2)}
\define\ot#1{\underset#1\to{\otimes}}
\define\ota{\ot{\g A}}
\define\pot#1{\underset#1\to{\widehat\otimes}}
\define\pott{\widehat\otimes}
\define\pota{\pot{\g A}}

\define\tens#1{#1\pot{#1}#1}

\define\Coh#1#2{\Cal H^1(#1,#2)}
\define\Ho#1#2{\Cal H^{#1}(#2,#2^*)}
\define\HoM#1{\Cal H^1_{\g M}(#1,#1^*)}

\define\lp#1{\ell_p(#1)}
\define\LP#1{L_p(\mu,#1)}

\define\clspan{\operatorname{clspan}}

\define\set#1{\{#1\}}

\topmatter
\title Push-outs of derivations\endtitle
\author Niels Gr\o nb\ae k\endauthor

\abstract
 Let $\g A$ be a Banach algebra and let $X$ be a Banach $\g
A$-bimodule. In studying $\Coh{\g A}X$ it is often useful to extend a
given derivation $D\colon\g A\to X$ to a Banach algebra $\g B$
containing $\g A$ as an ideal, thereby exploiting (or establishing)
hereditary properties. This is usually done using (bounded/unbounded)
approximate identities to obtain the extension as a limit of operators
$b\mapsto D(ba)-b.D(a)\,,\;a\in\g A$ in an appropriate operator
topology, the main point in the proof being to show that the limit map
is in fact a derivation. In this paper we make clear which part of
this approach is analytic and which algebraic by presenting an
algebraic scheme that gives derivations in all situations at the cost
of enlarging the module. We use our construction to give improvements
and shorter proofs of some results from the literature and to give a
necessary and sufficient condition that biprojectivity and biflatness
are inherited to ideals.
\endabstract

\endtopmatter
\document \heading0. Introduction\endheading

For a Banach algebra $\g A$ and a Banach $\g A$-bimodule $X$ the
computation of the bounded Hochschild cohomology group $\Coh{\g A}X$
can often be facilitated by extending a given derivation $D\colon\g
A\to X$ to some appropriate Banach algebra $\g B$ containing $\g A$ as
an ideal. Perhaps the first use of this is [J,Proposition 1.11] where
derivations from a Banach algebra with a bounded two-sided approximate
identity into a so-called neounital module is extended to the
multiplier algebra. In the case of group algebras $L_1(G)$ this gives
the link between dual Banach $L_1(G)$ modules and w$^*$-continuous
group actions which is essential in the proof that $L_1(G)$ has
trivial cohomology with coefficients in dual modules if and only if
the group $G$ is amenable [J, Theorem 2.5]. In [G1] extension
techniques are used to study weak amenability of commutative Banach
algebras. In [Theorem 7.1,GLo] so-called approximate and essential
amenability are established for certain abstract Segal algebras by
means of an extension result. In [GLa2] a similar argument is used on
symmetric Segal algebras on a SIN group to establish approximate weak
amenability and on an amenable group to establish approximate
amenability. In the study of weak amenability of Banach algebras $\Cal
A(X)$ of approximable operators on a Banach space $X$ [Bl, G2] it is
exploited that derivations can be extended to the Banach algebra $\Cal
B(X)$ of all bounded operators on $X$. All these extension techniques
are variations over the same extension problem.  The data of this
problem are a Banach algebra $\g A$ and a continuous injection as an
ideal into a larger Banach algebra $\g B$ and a bounded derivation
$D\colon \g A\to X$ into a Banach $\g A$-bimodule. The corresponding
diagram problem is
$$
\CD
\g A@>\iota>>\g B\\
@V D VV     @VV ? V\\
X@> ? >>?
\endCD\quad,\tag DP
$$
where we are looking for a Banach $\g B$-bimodule $\wt X$, a bounded $\g
A$-bimodule map $X@>\iota>>\widetilde
X$ and a bounded
derivation $\widetilde D\colon \g B\to\widetilde X$ such that the
completed diagram is commutative. Provided we can extend the module
action of $\g A$ on $X$ to an action of $\g B$, the problem (DP) is purely
algebraic (the boundedness of $\widetilde D$ will be obvious from the
construction) and can always be solved. Under various (mainly) topological
conditions we shall further investigate exactly how the modules $X$ and
$\widetilde X$ are related and how information about $\widetilde D$ is
transfered to $D$.

\heading{1. The Construction}\endheading

First we recall standard concepts, notation and elementary
results. For further details the reader is referred to [H]. Throughout
$\g A$ is a Banach algebra and $X$ and $Y$ are Banach $\g
A$-bimodules. The $\g A$-balanced projective tensor product is denoted
$-\pota-$. A Banach $\g A$ bimodule is called {\it induced} if
exterior multiplication induces a bimodule isomorphism $\g A\pota
X\pota\g A\mapsto X$. In particular $\g A$ is {\it self-induced} if
multiplication induces an isomorphism $\g A\pota \g A\mapsto \g A$. We
note that, if $\g A$ has a bounded approximate identity (BAI), then a
Banach $\g A$ module is {\it neounital} in the sense of [J], if and
only if it is induced. The vector space of bounded left (right) Banach
module homomorphisms $\g A\to X$ is denoted $\ah X$ ($\ha X$). These
are Banach $\g A$-bimodules with the uniform operator norm and module
operations given by

$$
\align
&a.S=a.S(\cdot),\;S.a=S(a\cdot),\quad a\in\g A,S\in \ha X\\
&a.T=T(\cdot a),\;T.a=T(\cdot).a,\quad a\in\g A,T\in \ah X
\endalign
$$ 

We note the important instances of hom-tensor duality 
$$
(\g A\pota X)^*\cong\ha{X^*},\;(X\pota\g A)^*\cong\ah{X^*}
$$ 
with isometric bimodule identifications thru
$$
<a\ota x,\Phi>=<x,S(a)>,\;<x\ota a,\Psi>=<x,T(a)>$$
for $ a\in \g A,x\in X,\Phi\in(\g A\pota X)^*,\Psi\in(X\pota\g
A)^*,S\in\ha{X^*},T\in\ah{X^*}.$

We define the {\it left} and {\it right annihilators} and the {\it
annihilator} of $X$ as
$$
\align
\annl X&=\{x\in X\mid \g A.x=\{0\}\}\\
\annr X&=\{x\in X\mid x.\g A=\{0\}\}\\
\ann X&=\annl X\cap\annr X\,.
\endalign
$$
 
A {\it module derivation} is a bounded linear map $D\colon\g A\to X$
satisfying the derivation rule
$$
D(ab)=a.D(b)+D(a).b,\quad a,b\in\g A.
$$
A derivation is {\it inner} if it is of the form $a\mapsto a.x-x.a$
for some $x\in X$.
The vector space of all derivations $D\colon\g A\to X$ is denoted
$\Cal Z^1(\g A,X)$ and the vector space of all inner derivations is
denoted $\Cal B^1(\g A,X)$. The {\it (first) Hochschild cohomology group} of
$\g A$ with coefficients in $X$ is $\Coh{\g A}X=\Cal Z^1(\g
A,X)\bigl/\Cal B^1(\g A,X)$.

We start by generalizing the notion of double centralizer ([J]) to
module homomorphisms, see also {\it multipliers} in [S]. 

\definition{Definition 1.1} 
 A {\it double centralizer} from $\g A$ into $X$ is a pair $(S,T)\in \ha
 X\times\ah X$ satisfying
$$
a.S(\alpha)=T(a).\alpha,\quad a,\alpha\in\g A\,.
$$
\enddefinition

\flushpar We note 

Each $x\in X$ gives rise to a double centralizer $(L_x,R_x)$, where
$L_x\colon a\mapsto x.a$ and $R_x\colon a\mapsto a.x$. 

For a double centralizer $(S,T)$ and $a\in \g A$ we have $a.S=L_{Ta},\,
a.T=R_{Ta},\, S.a=L_{Sa},\, T.a=R_{Sa}$.

Thus, norming $\ha X\times\ah X$ by $\norm{(U,V)}=\max\{\norm U,\norm
V\}$, the set of double centralizers is a Banach sub-bimodule of $\ha
X\times\ah X$, which we denote $\DC X$. The map $\iota_X\colon
x\mapsto (L_x,R_x),\; x\in X$ is a bounded bimodule homomorphism of
$X$ into $\DC X$ with the property that  $\g A.\DC X+\DC X.\g
A\sbst\iota_X(X).$
  
\medskip
With the definition to follow we have stipulated a setting that covers
most applications.

\definition{Definition 1.2} Let $(\g A,\vert\cdot\vert)$ and $(\g
B,\norm\cdot)$ be Banach algebras. We say that $\g B$ is an {\it
envelope} of $\g A$ if

\roster
\item $\g A$ is an ideal in $\g B$ with bounded inclusion $\iota\colon\g
A\to\g B$;
\item there is $C>0$ so that $\max\{\vert ba\vert,\vert ab\vert\}\leq C\vert a\vert\norm
b\,,\;a\in \g A, b\in\g B$.
\endroster
\enddefinition

We can now solve the diagram problem (DP). In cases of interest the
$\g A$-module action on $X$ can be extended to an action of $\g B$, so
we shall assume from the outset that $X$ is a $\g B$-bimodule.

\proclaim{Theorem 1.3} Let $\g B$ be an envelope of\/ $\g A$ and let $X$ be
a Banach $\g B$-bimodule. Then $X$ is naturally a Banach $\g
A$-bimodule and $\DC X$ is naturally a Banach $\g B$-bimodule. To each
bounded derivation $D\colon\g A\to X$ there is a  bounded
derivation $\widetilde D\colon\g B\to\DC X$, such that the diagram
$$
\CD
\g A@>\iota>>\g B\\
@V D VV     @VV\widetilde D V\\
X@>\iota_X>>\DC X
\endCD.
$$
is commutative.
Furthermore:

Each $(S,T)\in\DC X$ determines a derivation $a\mapsto
S(a)-T(a), \; a\in\g A$ of\/ $\g A$ into $X$. 

Suppose that $\overline{\g A^2}=\g A$. Then $\wt D$ is uniquely
determined by $D$ and in case $\wt D$ is inner,  there is
$(S,T)\in\DC X$ such that $D(a)=S(a)-T(a),\; a\in\g A$.
\endproclaim
\demo{Proof} We make $\DC X$ a Banach
$\g 
B$-bimodule by defining 
$$
b.(L,R)=(^bL,{}^bR)\text{ and }(L,R).b=(L^b,R^b),\quad b\in\g B,
(L,R)\in\DC X,
$$
with
$$
^bL(a)=b.L(a), {}^bR(a)=R(ab), L^b(a)=L(ba),R^b(a)=R(a).b,\quad a\in\g A\,.
$$
We define $\widetilde D$ by
$$
\widetilde D(b)=(\Cal L(b),\Cal R(b)),\quad b\in\g B\,,
$$
where
$$
\Cal L(b)(a)= D(ba)-b.D(a), \Cal R(b)(a)=D(ab)-D(a).b\;, a\in\g A,b\in\g
B\,.
$$
 We check that we can complete the diagram as wanted. Let
$a,c\in \g A$ and $b,b_1,b_2\in\g B$.
\flushpar The double centralizer property:
$$
\align
\Cal L(b)(ac)=D(bac)-b.D(ac)&=ba.D(c)+D(ba).c-b.(a.D(c)+D(a).c)\\
&=(D(ba)-b.D(a)).c=\Cal L(b)(a).c\;,
\endalign
$$
$$
\align
\Cal R(b)(ac)=D(acb)-D(ac).b&=a.D(cb)+D(a).cb-(D(a).c+a.D(c)).b\\
&=a.(D(cb)-D(c).b)=a.\Cal R(b)(c)\;,
\endalign
$$
and
$$
\align
&a.\Cal L(b)(c)=a.(D(bc)-b.D(c))=D(abc)-D(a).bc
-ab.D(c)=\\
&D(ab).c-D(a).bc=(D(ab)-D(a).b).c=\Cal R(b)(a).c
\endalign
$$
The derivation property:
$$
\align
(^{b_1}\Cal L(b_2)+\Cal L(b_1)^{b_2})(a)&=b_1.\Cal L(b_2)(a)+\Cal
L(b_1)(b_2a)\\
&=b_1.(D(b_2a)-b_2.D(a))+D(b_1b_2a)-b_1.D(b_2a)\\
&=D(b_1b_2a)-b_1b_2.D(a)\\
&=\Cal L(b_1b_2)(a)\;,
\endalign
$$
and
$$
\align
(^{b_1}\Cal R(b_2)+\Cal R(b_1)^{b_2})(a)&=\Cal R(b_2)(ab_1)+\Cal
R(b_1)(a).b_2\\
&=D(ab_1b_2)-D(ab_1).b_2+(D(ab_1)-D(a).b_1).b_2\\
&=D(ab_1b_2)-D(a).b_1b_2\\
&=\Cal R(b_1b_2)(a)\;,
\endalign
$$
Since clearly $(L_{D(a)},R_{D(a)})=(\Cal L(a),\Cal R(a))$ the diagram
is commutative.

\smallskip

It is immediate to verify that a each $(S,T)\in\DC X$ defines a
derivation by $a\mapsto S(a)-T(a)$, see also [GLa2].

Now suppose that $\overline{(\g A)^2}=\g A$, and let $\wt D\colon b\mapsto(\Cal
L(b),\Cal R(b))$ be any derivation that makes the diagram commutative. From the derivation
identities 
$$
L_{D(ab)}={}^a\Cal
L(b)+L_{D(a)},\;R_{D(bc)}={}^bR_{D(c)}+\Cal R(b)^c\;\;a,c\in\g A,b\in
\g B
$$
we get $\Cal L(b)(ac)=D(bac)-b.D(ac)$ and $\Cal
R(b)(ac)=D(acb)-D(ac).b$ so $\wt D$ is uniquely
determined by $D$.

If $\wt D(a)=(S,T).a-a.(S,T),\;a\in\g A$, then
$$
(L_{D(a)},R_{D(a)})=(L_{S(a)-T(a)},R_{S(a)-T(a)}),\;a\in\g A,
$$
so for all $a_1,a_2\in\g A$ we have
$$
\align
&D(a_1).a_2=(S(a_1)-T(a_1)).a_2=S(a_1a_2)-T(a_1).a_2;\\
&a_1.D(a_2)=a_1.(S(a_2)-T(a_2))=T(a_1).a_2-T(a_1a_2)
\endalign
$$
so that $D(a_1a_2)=S(a_1a_2)-T(a_1a_2)$. Since $\overline{\g A^2}=\g
A$, we get $D(a)=S(a)-T(a),\;a\in\g A$.
\enddemo

Very often one is interested in dual modules. The property of being a
dual module is preserved by the double centralizer construction. (This also follows
easily from the description as multipliers in [S].)

\proclaim{Proposition 1.4} Let $X^*$ be the dual of a Banach $\g A$-bimodule $X$.
Then there is an isometric bimodule isomorphism
$$
\DC{X^*}\cong((\g A\pota X\oplus X\pota \g A)/N)^*
$$ 
where $\oplus$ denotes $\ell_1$-direct sum and $N=\clspan\{(a\ota
x.\alpha,-a.x\ota\alpha)\mid a,\alpha\in \g A,\;x\in X\}$. If $q\colon
\g A\pota X\oplus X\pota \g A\to(\g A\pota X\oplus X\pota \g A)/N$ is
the quotient map and 
$\mu\colon(\g A\pota X\oplus X\pota \g A)\to X$ is the canonical map
induced by $(a\ot{\g A}x,x'\ot{\g A}a')\mapsto a.x+x'.a'$, then
$q^*\circ\iota_{X^*}=\mu^*$. In particular, $\iota_{X^*}$ is bounded below,
if and only if $\mu$ is surjective.
\endproclaim

\demo{Proof} We note that $N\subseteq\ker \mu$. The statement is then
a consequence of the natural isometric isomorphisms $\ha{X^*}\cong(\g
A\pota X)^*$ and $\ah{X^*}\cong (X\pota\g A)^*$ and the definition of
double centralizers.
\enddemo

In the important case of induced modules the description is particularly nice.

\proclaim{Corollary 1.5} Suppose that $X$ is an induced $\g A$-module.
Then $\iota_{X^*}\colon X^*\to\DC{X^*}$ is an isomorphism.
\endproclaim
\demo{Proof} In view of the proposition above we must prove
$\ker\mu\subseteq N$. Since $X$ is induced a generic element in
$\ker\mu$ has the form $k=(\sum_na_n\ota x_n.c_n,\sum_na_n'.x_n'\ota
c_n')$ with $\sum_na_n\ota x_n\ota c_n=-\sum_na_n'\ota x_n'\ota c_n'$
in $\g A\pota X\pota\g A$. It follows that $\sum_na_n'.x_n'\ota
c_n'=-\sum_na_n.x_n\ota c_n$ so that $k=\sum_n(a_n\ota
x_n.c_n,-a_n.x_n\ota c_n)$, which is clearly in $N$.

\enddemo

In the following we investigate further properties of the map
$\iota_X$. Clearly $\iota_X$ is injective, if and only if $\ann
X=\nul$. It turns out that in case of injectivity the pair
$(\iota_X,\DC X)$ is a universal element (in a certain category, only
to be implicitly specified) with the consequence that our solution to
the diagram problem is actually a push-out. Henceforth we shall only
work with annihilator-free modules. To which extent this is a
limitation is illustrated by the following remark.

\remark{Remark 1.6} For any left Banach $\g A$ module $X$ there is a
smallest closed submodule $\Cal N$ such that $\annl {(X/\Cal
N)}=\nul$. To see this, let $\Cal M$ be the set of closed submodules
$N\sbst X$ such that
$$
\forall x\in X\quad\g A.x\sbst N\implies x\in N\,.
$$
Put $\Cal N=\bigcap_{N\in\Cal M}N$. Then one checks that $\Cal N$ has
the desired property. The module $\Cal N$ may be constructed by a
transfinite procedure as follows. For any left submodule $N\sbst X$ we
set $N\colon\g A=\{x\in X\mid \g A.x\sbst N\}$. Let $N_1=\nul\colon\g
A$. For a successor ordinal $\beta=\alpha+1$ we set
$N_\beta=N_\alpha\colon \g A$ and for a limit ordinal $\gamma$ we set
$N_\gamma=\overline{\bigcup_{\alpha<\gamma}N_\alpha }$. Then $\Cal
N=\bigcup_{\card\alpha\leq \card X} N_\alpha$. Note that by a
cardinality argument any given
sequence in $\Cal N$ is contained in some $N_\gamma$, so $\Cal N$ is
automatically closed. Replacing $\g A$ by its opposite algebra
or by its enveloping algebra we get the same statements for right
modules and bimodules.
\endremark

\medskip
These elementary observations are part of the following description of
the double centralizer module as a universal object.

\proclaim{Proposition 1.7} Let $X$ be a Banach $\g A$-bimodule with
$\ann X=\nul$. The map $\iota_X\colon x\mapsto (L_x,R_x)$ is a
bounded $\g A$-bimodule monomorphism from $X$ into $\DC X$. The
pair $(\iota_X,\DC X)$ is a universal object in the sense: If $\wt X$
is a Banach $\g A$-bimodule and $j\colon X\to \wt X$ is a bounded $\g
A$-bimodule monomorphism with
$$
\g A.\wt X+\wt X.\g A\sbst j(X)\,,\tag$*$
$$
then there is a unique bounded $\g A$-bimodule homomorphism $\wt j\colon
\wt X\to \DC X$ such that $\iota_X=\wt j\circ j$.
\endproclaim
\demo{Proof} Since $\ann X=\nul$ the assignment $x\mapsto (L_x,R_x)$ defines
a monomorphism which, by the paragraphs after Definition 1.1,
satisfies ($*$). Let $j\colon X\to\wt X$ be a monomorphism satisfying 
($*$). By this property $L_{\wt x}(\g A)+R_{\wt x}(\g A)\sbst j(X)$ for
any $\wt x\in\wt X$, so we may define $S=j^{-1}\circ L_{\wt x}$ and
$T=j^{-1}\circ R_{\wt x}$. By the closed graph theorem $S$ and $T$ are
bounded linear maps, and since $j$ is a bimodule homomorphism, the
pair $(S,T)\in \DC X$. Define $\wt j(\wt x)=(S,T)$. It is clear that
$\iota_X=\wt j\circ j$, and yet an application of the closed graph
theorem tells that $\wt j$ is bounded. It remains to establish
uniqueness. Hence 
suppose that $\wt j_1\circ j=\wt j_2\circ j=\iota_X$. Let $\wt x\in
\wt X$ and set $(S_i,T_i)=\wt j_i(\wt x),\;i=1,2$. For any $a\in\g A$
we have $(L_{T_ia},R_{T_ia})=a.\wt j_i(\wt x)=\wt j_i(a.\wt
x)=\iota_X(j^{-1}(a.\wt x)),\;i=1,2$. Thus
$\iota_X(T_1a)=\iota_X(T_2a)$, and, since $\iota_X$ is a
monomorphism, $T_1=T_2$. By considering $\wt j_i(\wt x).a$ we get
$S_1=S_2$ as well.

\enddemo

\heading{2. Examples}\endheading
The examples to follow serve to illustrate a unified approach to
various results in the literature about inherited cohomology.

\subheading{2.1 Amenability}

Extension techniques were first used in [J], where derivations were
extended from Banach algebras with bounded approximate identities to
their double centralizer algebras. Recall that a Banach algebra is
{\it amenable} if $\Coh{\g A}{X^*}=\{0\}$ for all Banach $\g
A$-bimodules $X$.
 
\proclaim{Theorem 2.1 {\rm ([J, Proposition 5.1])}} Let $\g B$ be
amenable and let $\g A\sbst\g B$ be a closed two-sided ideal. Then $\g
A$ is amenable, if (and only if) $\g A$ has a BAI.
\endproclaim

\demo{Proof} The only if part is a standard fact about amenable Banach
algebras. As observed in [J], to prove amenability of $\g A$, we only
have to look at derivations into modules $X^*$ where $X$ is a
neounital $\g A$-module. The module action on such a module is
naturally extended to an action of $\g B$ for instance by noting that
$X\cong\g A\pota X\pota\g A$. By Corollary 1.5 we have $X^*=\DC{X^*}$,
so the conclusion follows immediately by amenability of $\g B$.
\enddemo 

\subheading{2.2 Weak amenability} A Banach algebra is {\it weakly
amenable} if $\Coh{\g A}{\g A^*}=\nul$. For commutative Banach
algebras this is equivalent to $\Coh{\g A}X=\nul$ for all symmetric
modules $X$. In this situation the
hereditary properties of weak amenability are simple.

\proclaim{Theorem 2.2 {\rm ([G, Corollary 1.3 ])}} Let $\g A$ be a
commutative Banach algebra and let $\g B$ be a commutative weakly
amenable envelope. Then $\g A$ is weakly amenable, if (and only if)
$\overline{\g A^2}=\g A$
\endproclaim

\demo{Proof} Let $D\colon \g A\to \g A^*$ be a derivation. Since $\g
B$ is commutative $\DC{\g A}$ is symmetric, so the push-out derivation
$\wt D$ is the zero map by weak amenability of $\g B$. Since
$\overline{\g A^2}=\g A$, the map $\iota_{\g A^*}$ is a monomorphism,
whence $D=0$.

As above the \lq only if' part is a standard fact about weakly
amenable Banach algebras.   
\enddemo

\subheading{2.3 Approximately inner derivations} The notion of
approximately inner deri\-vations was introduced in [GLo]. A
derivation $D\colon\g A\to X$ is {\it approximately inner} if $D$
belongs to the strong operator closure of $\Cal B^1(\g A,X)$,
i.e\. if there is a net $x_\gamma\in X$ such that $D(a)=\lim_\gamma
a.x_\gamma-x_\gamma.a,\;a\in\g A$. The Banach algebra $\g A$ is termed
{\it approximately amenable\/}, if for all Banach $\g A$-bimodules $X$
all derivations $D\colon\g A\to
X^*$ are approximately inner, and {\it approximately weakly amenable}
if all derivations $D\colon\g A\to\g A^*$ are approximately inner. We
have the generalization of [GLo, Corollary 2.3]

\proclaim{Theorem 2.3} Let $\g A$ have an approximately amenable envelope
$\g B$. If $X$ is $\g A$-induced, then every derivation $D\colon\g
A\to X^*$ is approximately inner. In particular, if $\g A$ has a BAI,
then $\g A$ is approximately amenable. 
\endproclaim
\demo{Proof} The first part is an immediate consequence of Corollary
1.5. If $\g A$ has a BAI, we only have to consider derivations into
duals of neounital modules [GLo, Proposition 2.5]. But neounital
modules over Banach algebras with a BAI are exactly the induced
modules.
\enddemo

\subheading{2.4 Abstract Segal algebras}  Let $\g S\sbst\g A$ be an
abstract symmetric Segal algebra of 
a Banach algebra $\g A$, i.e $\g A$ is an envelope of $\g S$ and $\g S$ is
dense in $\g A$, see [Bu]. From Theorem 2.2 it follows that, if $\g A$ is
commutative and weakly amenable, then $\g S$ is weakly amenable when
$\overline{\g S^2}=\g S$, thus strengthening  Theorem 3.3 of [GLa1].

In [GLa2] derivations from concrete symmetric Segal algebras are
investigated. A concrete symmetric Segal algebra has an AI which is a
BAI for $L^1(G)$. The results Theorem 3.1 and Remark 3.4 of this
reference follow from

\proclaim{Theorem 2.4 {\rm ([GLa2])}} Suppose that $\g A$ is amenable and
$X$ is a Banach $\g A$-bimodule. Then every derivation $D\colon \g
S\to X^*$ has the form $D(a)=S(a)-T(a),\;a\in\g S$. If $\g S$ has an
AI which is a BAI for $\g A$, then every derivation $D\colon\g S\to X$
is approximately inner and in particular $\g S$ is approximately
weakly amenable.
\endproclaim
\demo{Proof} Since $\overline{\g S^2}=\g S$ the first part follows
from Theorem 1.3 and Proposition 1.4. For the proof about being
approximately inner, we may as in [GLa2] reduce the task to showing
that, if $X$ is neounital over $\g A$, then every derivation $D\colon\g
S\to X^*$ is a $w^*$ limit of inner derivations. This reduction is
obtained by an appeal to the reduction  to neounital modules in [J]
combined with a standard application of Goldstine's Theorem and
Hahn-Banach's Theorem. Hence let $(e_\gamma)_\Gamma$ be an AI for $\g
S$. Since $\g A$ is amenable, there is $(S,T)\in\DC{X^*}$ so that $\wt
D(a)=a.(S,T)-(S,T).a,\;a\in\g S$. Since $X$ is neounital, we get
$D(a)=w^*-\lim_\gamma a.T(e_\gamma)-T(e_\gamma).a,\; a\in\g S$.  
\enddemo

\subheading{2.5 Biprojective and biflat Banach algebras}

Biprojective and biflat Banach algebras can be described in terms of
derivations, as demonstrated in [S]. Biprojectivity and -flatness are
nicely inherited 
by ideals.
\proclaim{Theorem 2.5} Let $\g A$ be a closed ideal of a biprojective
[biflat] Banach algebra $\g B$. Then  $\g A$ is biprojective [biflat]
if and only if $\g A$ is self-induced.
\endproclaim

\demo{Proof} The property of being self-induced is well-known to hold
for all biflat, and hence all biprojective Banach algebras. By
[S, Theorem 5.9] $\g A$ is biprojective if and only if $\Coh {\g A}{\DC
X}=\nul$ for all Banach $\g A$ bimodules $X$. In fact, to establish
biprojectivity it suffices to check this for $X$ being the kernel of
the multiplication $\g A\pott\g A\to\g A$. Clearly this $X$ is a Banach $\g
B$-bimodule, so we may apply Theorem 1.3. By [S, Lemma 5.2(iv)] the map $\iota_{\DC
X}\colon\DC X\to\DC{\DC X}$ is an isomorphism, since $\g A$ is
self-induced. Using that $\g B$ is biprojective it is possible to
prove directly that $\Coh{\g B}{\DC X}=\nul$ and therefore that
$\Coh{\g A}{\DC X}=\nul$. However, we prefer to use Theorem 1.3 once more,
now applied to the identity map $\g B\to\g B$ and $\DC X$ viewed as a
$\g B$ module. Invoking the universal
property of the push-out we conclude that every derivation $D\colon\g
A\to \DC X$ factors  as $D=j\circ \wt D\circ \iota$ where $\wt D\colon
\g B\to\DC{\DC X}_{\g B}$ and $\iota\colon \g A\to\g B$ is the
inclusion. The subscript indicates that double centralizers are with
respect to $\g B$. Since $\g B$ is biprojective, $\wt D$ is inner, thus
giving that $D$ is inner.

By applying above proof to the dual module $X^*$ we get the statement
about biflatness.
\enddemo


\widestnumber\key{GLa2}
\Refs

\ref\key Bl \by Blanco, A. \paper Weak amenability of $\Cal A(E)$ and the
geometry of $E$\jour J. London Math. Soc. (2)\vol 66\yr 2002\pages
721--740
\endref

\ref \key{Bu} \by Burnham, J. T. \paper Closed ideals in subalgebras of
Banach algebras. I.\jour Proc. Amer. Math. Soc. \vol 32
\yr1972\pages551--555
\endref

\ref\key{GLa1}\by Ghahramani, F.; Lau, A. T.-M. \paper Weak
amenability of certain classes of Banach algebras without bounded
approximate identities\jour Math. Proc. Cambridge Philos. Soc. \vol
133 \yr2002\pages357--371
\endref

\ref\key{GLa2}\bysame \paper
Approximate weak amenability, derivations and Arens regularity of
Segal algebras\jour Studia Math. \vol 169 \yr2005 no. 2\pages 189--205
\endref

\ref\key{GLo}\by Ghahramani, F.; Loy, R. J. \paper Generalized
notions of amenability\jour  J. Funct. Anal.  \vol 208  \yr
2004)\pages 229--260
\endref

\ref\key{G1} \by Gr\o nb\ae k, N. \paper Commutative Banach algebras,
module derivations, and semigroups\jour J. London
Math. Soc. (2)\vol40\yr1989\pages 137--157
\endref

\ref\key{G2}\bysame\paper Self-induced Banach
algebras\jour Contemp. Math. \vol 263\yr 2004\pages129--143\endref

\ref\key{H}\by Helemski\u \i, A. Ya.\book  Banach and
Locally Convex Algebras\publ Oxford University Press\publaddr New
York\yr 1993\endref 

\ref \key{J}\by Johnson, B.E. \paper Cohomology in Banach
algebras\jour Mem. Amer. Math. Soc.\vol 127\yr 1972
\endref 

\ref\key{S}\by Selivanov, Yu. V.\paper Cohomological
Characterizations of Biprojective and Biflat Banach Algebras\jour
Mh. Math.\vol 128\yr 1999\pages 35--60
\endref

\endRefs
\enddocument